\documentclass{article}
\usepackage{amsmath,amstext,amsthm,amsfonts}
\usepackage[dvips]{graphicx}
\usepackage{epic,eepic}

\theoremstyle{plain}
\newtheorem{theorem}{Theorem }[section]
\newtheorem{proposition}[theorem]{Proposition}
\newtheorem{lemma}[theorem]{Lemma}
\newtheorem{corollary}[theorem]{Corollary}
\newtheorem{maintheorem}{Theorem}

\theoremstyle{definition}
\theoremstyle{remark}
\newtheorem{remark}[theorem]{Remark}

\newcommand{\field}[1]{\mathbb{#1}}
\newcommand{\real}{\field{R}}

\newcommand{\al} {\alpha}       
        
\newcommand{\ga} {\gamma}    
\newcommand{\de} {\delta}       
\newcommand{\ep} {\epsilon}

\newcommand{\la} {\lambda}      \newcommand{\La}{\Lambda}

\newcommand{\vfi}{\varphi}
\newcommand{\om} {\omega}

\newcommand{\mundo}{\operatorname{Diff}^1_{\omega}(M)}
\newcommand{\submundo}{\operatorname{Diff}^{1+\al}_{\omega}(M)}

\newcommand{\Mundo}{\mathfrak{X}^{1}_m(M)}
\newcommand{\Submundo}{\mathfrak{X}^{\infty}_m(M)}
\newcommand{\Subsubmundo}{\mathfrak{X}^{1+\al}_m(M)}
\newcommand{\Newmundo}{\mathfrak{X}^{k+\al}_{\omega}(M)}
\newcommand{\com}{\mathbb{P}_c(M)}

\newcommand{\N}{\mathbb{N}}
\newcommand{\R}{\mathbb{R}}
\newcommand{\eps}{\varepsilon}
\newcommand{\Lp}{\Lambda^+}
\newcommand{\Lm}{\Lambda^-}
\newcommand{\Sun}{W^u_g(p)}

\newcommand{\ov}{\overline}

\newcommand{\SF}{{\cal F}}

\newcommand{\SK}{{\cal K}}

\newcommand{\SO}{{\cal O}}
\newcommand{\SP}{{\cal P}}

\newcommand{\SR}{{\cal R}}
\renewcommand{\SS}{{\cal S}}

\newcommand{\SU}{{\cal U}}
\newcommand{\SV}{{\cal V}}

\begin{document}

\title{\Large{\textbf{A Pasting Lemma I: the case of vector fields}}}

\author{ Alexander Arbieto\footnote{A. Arbieto is supported by CNPq-Brazil and
Faperj-Brazil} 
\text{ and} Carlos Matheus\footnote{C. Matheus is supported by Faperj-Brazil} }

\date{June 20, 2003}

\maketitle

\begin{abstract} 
We prove a perturbation (pasting) lemma for conservative (and symplectic) 
systems. This allows us to prove that $C^{\infty}$ volume preserving vector
fields are $C^1$-dense in
$C^{1}$ volume preserving vector fields\footnote{After the conclusion of this
work, Ali Tahzibi pointed out to us that this result was proved
in 1979 by Carlos Zuppa, although his proof is different from ours.}. 
Moreover, we obtain that $C^1$ robustly transitive conservative
flows in three-dimensional manifolds are Anosov and we conclude that
there are no geometrical Lorenz-like sets for conservative flows. Also,
by-product of the version of our pasting lemma for conservative diffeomorphisms,
we show that $C^1$-robustly transitive conservative $C^2$-diffeomorphisms admits
a dominated splitting, thus solving a question posed by Bonatti-Diaz-Pujals. In
particular, stably ergodic diffeomorphisms admits a dominated splitting. 
\end{abstract}

\section{Introduction}
\label{intro}

This is a preliminary version of a work in progress. So, we have to explain why
we are writing it. The explanation is somewhat simple: in the present work we
proved some results concerning some ``pasting'' (perturbation) lemma which 
were used by another researchers, although there is not even a \emph{preprint}
version of the theorems. However, the disadvantage of this preliminary version
is that it is too close to a draft of the theorems that we hope to prove in a
(if possible) near by future. But, we point out that the results below are 
correct, although some of them need a refinement. In particular, the suggestions
of the readers of this ``pre-preprint'' certainly will help in order to improve
this material. 

The main technical result in this work is a general tool for perturbing $C^2$
conservative dynamical systems over (large) compact subsets of phase space,
obtaining a globally defined and still conservatove perturbations of the
original system. In fact, we give several versions of this pasting lemmas, for
discrete as well as for continuous time systems (and in this case is enough
suppose that the system is $C^1$). Before hand, let us describe some of its
consequences in dynamics.

\smallskip
{\it Robustly transitive diffeomorphisms}

A dynamical system is \emph{robustly transitive} if any $C^1$ nearby one has
orbits that are dense in the whole ambient space. Recall that a \emph{dominated splitting} is a continuous decomposition
$TM=E\oplus F$ of the tangent bundle into continuous subbundles which are
invariant under the derivative $Df$ and such that $Df|_E$ is more expanding/less
contracting than $Df|_F$ by a definite factor (see precise definitions and more
history on section~\ref{s.dom}). 

We prove: ``\emph{Let $f:M\to M$ a $C^2$ volume preserving diffeomorphism
robustly transitive among volume preserving diffeomoorphisms. Then $f$ amits a
dominated splitting of the tangent bundle}''. This is an extension of a theorem
by Bonatti, Diaz and Pujals.

\smallskip

{\it Robustly transitive vector fields}

First we prove an extension of a theorem by Doering in the conservative setting:
``\emph{Let $X$ be a divergent free vector field robustly transitive among
divergent free vector fields in a three dimensional manifold then $X$ is Anosov.}

It has been show by Morales-Pac\'{\i}fico and Pujals~\cite{MPP} that robustly
transitive sets $\Lambda$ of a 3-dimensional flows are \emph{singular hyperbolic}: tehre
exists a dominated splitting of the tangent bundle restricted to $\Lambda$ such
that the one-dimensional subbundle is hyperbolic(contracting or expanding and
the complementary one is volume hyperbolic. Robust transitivity means that
$\Lambda$ is the maximal invariant set in a neighborhood $U$ and, for any $C^1$
neasrby vector field, the maximal invariant set inside contains dense orbits. 

For generic dissipative flows robust transitivity can not imply hyperbolicity,
in view of the phenomenom of Lorenz-like attractors containing both equilibria
and regular orbits. However we also prove that Lorenz-like sets do not exist for
consevative flows: ``\emph{If $\Lambda$ is a robust transitive set within
divergent free vector fields then $\Lambda$ contains no equilibrium points.}

\section{The Pasting Lemmas}
\label{lemmas}
Before stating the main results of this section, the pasting lemmas, we recall
the following results from PDE's by Dacorogna-Moser~\cite[Theorem 2]{DM} on the 
divergence equation:

\begin{theorem}
\label{t.2dm}
Let $\Omega$ a manifold with $C^{\infty}$ boundary. Let $g\in
C^{k+\alpha}(\Omega)$ such that $\int_{\Omega}g=0$. Then there exist $v$ a
$C^{k+1+\alpha}$ vector field (with the same regularity at the boundary) such
that:
\begin{displaymath}
\left\{ \begin{array}{ll}
\text{div }v(x)=g(x) \text{, }x\in\Omega \\
v(x)=0 \text{, }x\in \partial \Omega \end{array}\right.
\end{displaymath}
Furthermore there exist $C=C(\alpha,k,\Omega)>0$ such that $\|v\|_{k+1+\alpha}\leq
C\|g\|_{k+\alpha}$. Also if $g$ is $C^{\infty}$ then $v$ is $C^{\infty}$.
\end{theorem}

\begin{remark}
In fact Dacorogna-Moser shows this result for open sets in $\R^n$. But the same
methods works for manifolds with boundaries, because it follows from the
solvability of $\triangle u=f$ with Neumann's condition (see~\cite[p.3]{DM}
or~\cite[p.265]{H}; see
also~\cite[Ch 4,Theorem 4.8]{Au},~\cite{GT} and~\cite{LU}).
\end{remark}

A similar result holds for diffeomorphisms, as showed by
Dacorogna-Moser~\cite[Theorem 1, Lemma 4]{DM} (see also~\cite{RY}):

\begin{theorem}\label{l.DM4}Let $\Omega$ a manifold with $C^{\infty}$ boundary.
Let $f,g\in
C^{k+\alpha}(\Omega)$ such that $f,g>0$. Then there exist $\vfi$ a
$C^{k+1+\alpha}$ diffeomorphism (with the same regularity at the boundary) such
that:
\begin{displaymath}
\left\{ \begin{array}{ll}
g(\vfi(x))\det(D\vfi(x))=\lambda f(x) \text{, }x\in\Omega \\
\vfi(x)=x \text{, }x\in \partial \Omega \end{array}\right.
\end{displaymath}
where $\lambda=\int g/\int f$.
Furthermore there exist $C=C(\alpha,k,\Omega)>0$ such that $\|\vfi\|_{k+1+\alpha}\leq
C\|f-g\|_{k+\alpha}$. Also if $f,g$ are $C^{\infty}$ then $\vfi$ is $C^{\infty}$.
\end{theorem}

\begin{remark}
A consequence of the Dacorogna-Moser's theorem is the fact of $\submundo$ is 
path connected. 
\end{remark}

We now state the pasting lemma:

\begin{theorem}[$C^1$-Pasting Lemma for vector fields]
\label{t.c1glue}
Given $\eps>0$ there exists $\delta>0$ such that if $X\in \Mundo$, $K$ a
compact set and $Y\in \Subsubmundo$ is $\delta$-$C^1$-close to $X$ on a 
small neighborhood $U$ of $K$ then there exist a $Z\in \Subsubmundo$ and $V$ such
that $K\subset V\subset U$ satisfying $Z|_V=Y$ and
$Z$ is $\eps$-$C^1$-close to $X$. If $Y\in\Submundo$ then $Z$ is also in
$\Submundo$.
\end{theorem}

\begin{proof}
Let $V$ a neighborhood of $K$ with $C^{\infty}$ boundary, compactly contained in
$U$ such that $V$ and $int(U^c)$ is a covering of $M$. Let $\xi_1$ and $\xi_2$ a
partition of unity subordinated to this covering such that $\xi_1|_V=1$ and
there exist $W\subset U^c$ such that $\xi_2|_W=1$. Let $\Omega=M\backslash
V\cup W$ a manifold with $C^{\infty}$ boundary and we fix $C$ as the constant
given by theorem~\ref{t.2dm}.

Now we choose $\delta_0$ such that if $Y_0$ is $\delta_0$-$C^1$-close to $X_0$ then
$T=\xi_1Y_0+\xi_2X_0$ is $\frac{\eps_0}{2}$-$C^1$-close to $X_0$ and if we get 
$g=\text{div }T$ then $g$ is $\frac{\eps_0}{2C}$-$C^1$-close to 0. Indeed, we note
that $g=\nabla \xi_1\cdot Y_0+\nabla \xi_2\cdot X_0$ is $C^1$-close to 0 if we get
$Y_0$ sufficiently $C^1$-close to $X_0$. Also we note that by the divergence theorem
(see the proof of the theorem~\ref{l.dense} below) $\int_{\Omega}g=0$. Clearly, $g$ is a 
$C^{\infty}$ function.

So we get $v$ the vector field given by the theorem~\ref{t.2dm} and extend as 0 in the rest of
the manifold this is a
$C^{\infty}$-vector field (because we have regularity of $v$ at the boundary)
and let $Z_0=T-v$ so this a $C^{\infty}$-vector field $\eps_0$-$C^1$-close to
$X_0$ which
satisfies the statement of the theorem.
\end{proof}

Using the same arguments of Dacorogna-Moser, we can produce perturbation in
higher topologies, provided that the original systems are smooth enough. We can sumarize
this as:

\begin{theorem}[Pasting Lemma with Higher Differentiability]
\label{l.gluinglemma}
Let $X\in\Newmundo$ ($k$ integer,
$0\leq\al
<1$ and $k+\al>1$). Fix $K\subset U$ a 
neighborhood of a compact set $K$. 
Let $Y$ be a vector field defined in $U_1\subset U$, $U_1$ an open set
containing $K$. If $Y\in\Newmundo$ is sufficiently $C^{r}$-close to $f$ for
$1\leq r\leq k+\al$ ($r$ is real)
then there exists some $C^{r}$-small 
$C^{k+\al}$-perturbation $Z$ of $X$ and an open set
$V\subset U$ containing $K$ such that $Z=X$ outside $U$ and $Z=Y$ on $V$.

Furthermore, there exists a constant $C=C(f,\dim(M),K,U,r)$, $\de_0$ and 
$\ep_0$ such 
that if $Y$ is $(C\cdot\de^{\widetilde{r}}\cdot\ep)$-$C^{r}$ close to $X$, 
$\de<\de_0$, $\ep<\ep_0$ then 
$Z$ is $\ep$-$C^{r}$ close to $X$ and $\textrm{support}(Z-X)\subset
B_{\de}(K)$ (= the $\de$-neighborhood of $K$). 
\end{theorem}

An interesting version of this result is the following corollary:

\begin{lemma}[The $C^{1+\alpha}$-Pasting Lemma for Vector Fields]\label{t.c2glue}
Given $\eps_0>0$ there exists $\delta_0>0$ such that if $X_0\in \Subsubmundo$, $K$ a
compact set and $Y_0\in \Subsubmundo$ is $\delta_0$-$C^1$-close to $X_0$ on a 
small neighborhood $U$ of $K$ then there exist a $Z_0\in \Subsubmundo$, and $V$ and
$W$ such
that $K\subset V\subset U\subset W$ satisfying $Z_0|_V=Y_0$, $Z_0|_{int(W^c)}=X_0$ and
$Z_0$ is $\eps_0$-$C^1$-close to $X_0$. Furthermore, if $X_0$ and $Y_0$ are
$C^{\infty}$ then $Z_0$ is also $C^{\infty}$.
\end{lemma}

\subsection{A Weak Pasting Lemma in the Conservative Case}

In this subsection we show a weak version of the pasting lemma for $C^2$
conservative systems, and then we apply it to the study of dominated splittings
for conservative systems.

The pasting lemma in this case is:

\begin{lemma}\label{l.C2pasting} If $f$ is a $C^2$-conservative diffeomorphism
and $x$ is a point
in $M$, then for any $0<\alpha <1,\eps >0$, there exists a
$\eps$-$C^{1}$-perturbation $g$ (which is a $C^{1+\alpha}$ diffeomorphism) of
$f$ such that, for some small neighborhoods $U\supset V$ of $x$, $g|_{U^c}=f$
and $g|_V=Df(x)$. 
\end{lemma}

\begin{proof} First consider a perturbation $h$ of $f$ such that
$h(y)=(1-\rho(y))(f(x)+Df(x)(y-x))+\rho(y) f(y)$ (in local charts), where $\rho$
is a bump function such that $\rho|_{B_{r/2}(x)}=1, \ \rho|_{M-B_r(x)}=0,
|\nabla\rho|\leq C/r\text{ and }|\nabla^2\rho|\leq C/r^2$.
Denote by $g(y)=f(x)+Df(x)(y-x)$. Now, we note that $||h-f||_{C^1}\leq C\cdot
||f||_{C^2}\cdot r$ and $||h-f||_{C^2}\leq C\cdot||f||_{C^2}$, where $C$ is a
constant. Then if we denote by $\theta(y)=\det Dh(y)$ the density function of
$h$, by the previous calculation, we obtain that $\theta$ is $C\cdot
||f||_{C^2}\cdot r$-$C^0$-close to $1$ (the density function of $f$) and
$||\theta -1||_{C^1}\leq C\cdot ||f||_{C^2}$. By the classical property of
convexity of the H\"older norms, we obtain that $||\theta -1||_{\alpha}\leq C
||f||_{C^2}\cdot r^{1-\alpha}$. So, applying the lemma~\ref{l.DM4}, we obtain
that there exists a $C^{1+\alpha}$ diffeomorphism $\chi$ which is $C\cdot
||f||_{C^2}\cdot r^{1-\alpha}$-close to the identity. To conclude the proof is
suffucient to note that $h\circ\chi$ is conservative and is close to $f$
provided that $r$ is small.
\end{proof}

\begin{remark} In the last section of this paper, we apply this weak pasting
lemma to prove that robust transitive conservative diffeomorphism admits a
dominated splitting. We hope that this lemma can be applied to obtain some other
interesting consequence. In fact, recently, Bochi-Fayad-Pujals were able to show
that stably ergodic systems are (generically) non-uniformly hyperbolic
(see~\cite{BFP}).  
\end{remark}

\subsection{A Weak Pasting Lemma in the Symplectic Case}

The above perturbation lemma can be done for symplectic diffeomorphisms without the
lost of the structure outside the neighboorhood of a periodic orbit, as
follows. Also we have the same result as in the conservative vector fileds case,
if we require more differentiability of the diffeomorphism .

\begin{lemma}
If $f$ is a $C^k$-symplectic
diffeomorphism ($k\geq 1$), $x\in Per^n(f)$ is a periodic point of $f$ and $g$
is a local diffeomorphism ($C^k$-close to $f$) defined in a small 
neighborhood $U$ of
the $f$-orbit of $x$,
then there exists a $C^k$-symplectic diffeomorphism $h$ ($C^k$-close to $f$) 
and some neighborhood $U\subset V$ of $x$ satisfying $h|_U = g$ and
$g|_{V^c}=f$
\end{lemma}

\begin{proof} Consider the $f$-orbit $\SO(x)$ of $x$. Since all the perturbations are
local, by Darboux's theorem, we can use local coordinates near
each point in $\SO(x)$, say $\vfi_i:U_i\rightarrow V_i$, where $f^i(x)\in U_i$
and $0\in V_i = B_{\ga_i}(0)\subset\real^n$. With respect to this coordinates,
we have the local maps $f:B_{\de_i}(0)\rightarrow V_{i+1}$, where $V_n=V_1$,
$B_{\de_i}(0)\subset V_i$ are small neighborhoods of $0$, $f(0)=0$. In the
symplectic case, we fix some bump function $\la:\real\rightarrow [0,1]$ such
that $\la(z)=1$ for $z\leq 1/2$ and $\la(z)=0$ for $z\geq 1$. Let $S_{1,i}$ a
\emph{generating function} for $g:B_{\de_i}(0)\rightarrow V_{i+1}$ and
$S_{0,i}$ a generating function for $f:B_{\de_i}(0)\rightarrow V_{i+1}$
(see~\cite{X} for definitions and properties of generating functions). Finally,
let $S(x,y)=\la(\frac{2||(x,y)||}{\de_i})\cdot S_{1,i}(x,y) + \big[
1-\la(\frac{2||(x,y)||}{\de_i}) \big]\cdot S_{0,i}(x,y)$. In particular, if $h$
the symplectic map associated to $S$, then $h=g$ in $B_{\de_i/4}(0)$
and $g=f$ outside $B_{\de_i/2}(0)$. To summarize, this construction give us a
symplectic $\ep$-$C^1$-perturbation $g$ of $f$, if the numbers $\de_i$
are small, such that $h$ coincides with $g$,
near to each point $f^i(x)$. 
\end{proof}

\begin{remark} We can not obtain global versions of the pasting lemma for
symplectic diffeomorphism, as proved for conservative ones, since we used a
local tool: \emph{generating functions}. 
\end{remark}

\section{Denseness Results}
\label{perturbations}

We start this section proving that $C^{\infty}$ conservative vector fields are $C^1$-dense in
$C^1$ conservative vector fields:

\begin{theorem}
\label{l.dense}
$\Submundo$ is $C^1$-dense on $\Mundo$.
\end{theorem}

\begin{proof}
Locally we proceed as following: take conservative local charts (Moser's
theorem) and get $\eta_{\eps}$ a $C^{\infty}$ Friedrich's mollifier. If
$X=(X_1,\dots, X_n)$, we define $X_{\eps}=(X_1*\eta_{\eps},\dots
X_n*\eta_{\eps})$ a $C^{\infty}$ vector field. Now we note the
$\frac{d}{dx_i}(X_i*\eta_{\eps})=(\frac{d}{dx_i}X_i)*\eta_{\eps}$, so $\text{div
}X_{\eps}=0$. Furthermore $X_{\eps}$ converges to $X$ in the $C^1$-topology. 

\smallskip

Now we get local charts $(U_i,\phi_i)$, $i=1,..., m$ as above and $\xi_i$ a
partition of unity subordinated to $U_i$, and open sets $W_i$ such that
$W_i\subset V_i:=supp(\xi_i)\subset U_i$
satisfies $\phi_i(W_i)=B(0,1/3)$ and $\phi_i(V_i)=B(0,2/3)$,
$\xi_i|_{W_i}=1$ and $\Omega:=M\backslash int(W_i)$ is a manifold with boundary
$C^{\infty}$. And we fix $C$ as the constant given by theorem~\ref{t.2dm}

In each of those charts 
we get $X_i$ a $C^{\infty}$ conservative vector fields $C^1$-close to $X$ such that if 
we take $Y=\sum_i \xi_iX_i$ a $C^{\infty}$ vector field and we denote $g=\text{div
 }Y$ a $C^{\infty}$ function then:
 
\begin{itemize} 
\item $g$ is $\frac{\eps}{2C}$-$C^1$-close to
0 (where $K$ is the constant on the above theorem). Indeed, we have:

$$g=\sum_i \nabla \xi_i \cdot X_i + \xi_i\cdot \text{div }X_i=\sum_i \nabla \xi_i
\cdot X_i.$$

So, it's sufficient take the $X_i$'s $C^1$-close enough to $X$.
\item $Y$ is $\frac{\eps}{2}$-$C^1$-close to $X$.

\end{itemize}

By the divergence theorem we get that $\int_{\Omega}g =0$, since:

$$\int_{\Omega} \text{div }Y=\int_{\partial \Omega} Y\cdot N=-\sum_i
\int_{\partial W_i} X_i\cdot N =-\sum_i \int_{W_i} \text{div } X_i=0.$$

Now we get $v$ the vector field given by theorem~\ref{t.2dm} and we define
$Z=Y-v$. We observe that $Z$ is a $C^{\infty}$ vector field because $v$ is
$C^{\infty}$ at the boundary. Also $\text{div }Z=0$. Finally since $Y$ is
$\frac{\eps}{2}$-$C^1$-close to $X$ and $v$ is $\frac{\eps}{2}$-$C^1$-close to 0
we get $Z$ $\eps$-$C^1$-close to $X$. 
The proof is complete.
\end{proof}

\begin{remark}As Ali Tahzibi pointed out to the second author after the
conclusion of this work, the previous 
theorem was proved by Zuppa~\cite{Zu} (in 1979) with a different proof. In fact,
since Dacorogna-Moser's result was not available at that time, Zuppa uses that
the Laplacian operator admits a right inverse. 
\end{remark}
In the next sections we study some other consequences of the pasting lemmas for
conservative systems.

\section{Robustly Transitive Conservative Flows in 3-manifolds are Anosov}
\label{mpp}
Let $M^3$ be a compact three manifold and $\om$ a smooth volume in $M$. A
vector field $X\in \Mundo$ is $C^1$-\emph{robustly transitive} (in 
$\Mundo$) if there is $\ep>0$ such that every $Y\in\Mundo$ 
an $\ep$-$C^1$-perturbation of $X$ is transitive.

In this section we study conservative flows in 3-dimensional manifolds. e
Actually, we have the following 
theorem by Morales, Pac\'{\i}fico and Pujals~\cite{MPP}: {\it Any isolated singular
$C^1$ robustly transitive set is a proper attractor.} This result is related
with a theorem by C. Doering~\cite{D}: {\it A $C^1$ robustly transitive flow is
Anosov}. In the same spirit of the diffeomorphism case, we can use the
pertubation lemma to prove the conservative version of these results.

More precisely, we prove that: \emph{$C^1$ robustly transitive conservative
$C^1$-flows in three-dimensional manifolds are Anosov}. With this theorem and
using the same techniques we obtain: \emph{If $\Lambda$ is an isolated robustly transitive set of a conservative
vector field $X\in \Mundo$ then it cannot have a singularity.} Singular robustly
transitive flows in 3-manifolds are related with geometrical Lorenz-like sets,
so our results gives that: \emph{There are no conservative Lorenz-like sets}. We
now state the theorems:

\begin{maintheorem}
\label{t.doering}
Let $X\in \Mundo$ be a conservative $C^1$-robustly transitive vector field in a
3-dimensional compact manifold. Then $X$ is Anosov.
\end{maintheorem}

\begin{maintheorem}
\label{t.e}
If $\Lambda$ is an isolated robustly transitive set of a conservative
vector field $X\in \Mundo$ then it cannot have a singularity.
\end{maintheorem}

In particular, because the geometrical Lorenz sets are robustly transitive and
carrying singularities~\cite{Max}, as a corollary we have:

\begin{corollary}
There are no geometrical Lorenz sets in the conservative setting.
\end{corollary}

\subsection{Some Lemmas}

As in the diffeomorphism case we cannot have elliptic periodic orbits (or
singularities) in the presence of robust transitivity. More precisely we have
the following lemmas:

\begin{lemma}
\label{l.noeliptic}
If $X\in \Mundo$ is $C^1$-robustly transitive then there are no elliptic
singularities (i.e. the spectrum of $DX(p)$ does not contain in $S^1$).
\end{lemma}

\begin{proof}
Suppose that there exist an elliptic singularity $\sigma$. Then we get $K$ a
small neighborhood of $\sigma$ such that $X_K$ is $C^1$-close to the linear
vector field $DX(\sigma)$. Now, we use the theorem~\ref{t.c1glue} to obtain a
vector field $Y$ $C^1$-close to $X$ with a (possibly smaller) invariant
neighborhood. This contradicts transitivity.
\end{proof}

We also have the same result for periodics orbits:

\begin{lemma}
\label{l.noelipticvector}
If $X\in \Mundo$ is $C^1$-robustly transitive then every periodic orbit is
hyperbolic.
\end{lemma}

\begin{proof}
Suppose that there exist an elliptic periodic orbit $p$. Then we get $K$ a
small neighborhood of $O(p)$ (e.g., a tubular neighborhood) such that
the linear flow induced by the periodic orbit restricted to $K$ is $C^1$-close
to $X$. 
Now, we use the theorem~\ref{t.c1glue} to obtain a 
vector field $Y$ $C^1$-close to $X$ with a (possibly smaller) invariant
neighborhood (an invariant torus). This contradicts transitivity.
\end{proof}

\begin{remark}
\label{r.lemapraconjunto}
The above corollaries holds for robustly transitive sets instead for the whole
manifold with only minor modifications on the statement.
\end{remark}

\subsection{Proof of the Theorem~\ref{t.doering}}

By lemma~\ref{l.noeliptic} there are no elliptic singularities. By an argument
of Robinson~\cite{Rob} we can supose that any singularity is hyperbolic.
Also by lemma~\ref{l.noelipticvector} every periodic orbit is hyperbolic.

\begin{remark} This needs
the robust transitivity and the perturbation lemma since Kupka-Smale's theorem
is false on dimension three 
for conservative vector fields.
\end{remark}

Now we note that the result below by Doering~\cite{D} also holds in the conservative
case. 
\begin{theorem}
\label{t.doeringcoisado}
Let $\SS(M)=\{X\in\Mundo; \text{ every }\sigma\in Crit(X) \text{ is
hyperbolic}\}$. If $X\in\SS(M)$ then $X$ is Anosov.
\end{theorem}

Now if $X$ is robustly transitive then we already know that any critical element
is hyperbolic so we have the result.

We only sketch the Doering's proof, since is straighforward that it works in the 
conservative case:
\begin{lemma}
\label{l.doe}
If $X\in\SS(M)$ then $M-Sing(X)$ has dominated splitting. In particular there
are no singularities.
\end{lemma}

\begin{proof} 

We start with a result on~\cite[Proposition 3.5]{D}: 

{\bf Claim:} \textit{There exist $c>0$ and $0<\lambda<1$ and $\SU$ neighborhood of $X$ in
$\SS(M)$ such that there exist $(C,\lambda)$-dominated splitting for any $q\in
Per(Y)$} (for the definition of $(C,\lambda)$ dominated splitting see~\cite{D}).

Let $x$ a regular point for any $Y\in\SU$, so there are $Y_n$ such that $x\in
Per(Y_n)$ (Pugh's Closing Lemma). Now we get the $P^{Y_n}$-hyperbolic splitting
over $O_{Y_n}(x)$ and by compacity of the Grassmanian we have a splitting over
$O_X(x)$ (by saturation) then by the Claim we have that this splitting is a
$(C,\lambda)$-dominated splitting. Now by the abstract lemma of
Doering~\cite[Lemma 3.6]{D} we obtain a dominated splitting on $M-Sing(x)$.

By hyperbolicity $x_0$ is an interior point of $(M-Sing(X))\cup x_0$, but this
is a contradiction with the domination by the abstract~\cite[Theorem 2.1]{D}.
\end{proof}

The proof finish with the following result that is also available in the
conservative setting.

\begin{lemma}[Liao's theorem]
\label{l.liao}
If $X\in\SS(M)$ without singularities then $X$ is Anosov.
\end{lemma}

\subsection{Proof of the Theorem~\ref{t.e}}

We follow the steps of~\cite{MPP}. We denote by $\SU$ a neighborhood such that
every $Y\in\SU$ is transitive.

As in the proof of the previous theorem, we can supose that any singularity or
periodic orbit is
hyperbolic. The first step is to prove that in fact it is a Lorenz-like
singularity. Recall that a singularity $\sigma$ is \emph{Lorenz-like} if its eigenvalues
$\lambda_2(\sigma)\leq\lambda_3(\sigma)\leq\lambda_1(\sigma)$ satisfies either $\lambda_3(\sigma)<0$ $\Rightarrow$
$-\lambda_3(\sigma)<\lambda_1(\sigma)$ or $\lambda_3(\sigma)>0$ $\Rightarrow$
$-\lambda_3(\sigma)>\lambda_2(\sigma)$.

\begin{lemma}
\label{l.lorenzlike}
If $X\in \Mundo$ is $C^1$-robustly transitive vector field then any
singularity $\sigma$ is Lorenz-like. 
\end{lemma}

\begin{proof}
First of all, the eigenvalues of $\sigma$ are real. Indeed, if $\omega=a+ib$ is an
eigenvalue of $\sigma$, then the others are $a-ib$ and $-2a$. Since the
singularity is hyperbolic, we can suppose also that $X\in \Submundo$ (using
theorem~\ref{l.dense}). By connecting 
lemma~\cite{WX}, we can assume that there exist a loop $\Gamma$ associated to
$\sigma$ which is a Shilnikov Bifurcation. Then by~\cite[page 338]{BS} there
is a vector field $C^1$-close to $X$ with an elliptic singularity which gives a
contradiction (We observe
that the regularity for the bifurcation in~\cite{BS} is more than or
equal to seven).

Let $\lambda_2(\sigma)\leq\lambda_3(\sigma)\leq\lambda_1(\sigma)$ the
eigenvalues. Now, $\lambda_2(\sigma)<0$ and $\lambda_1(\sigma)>0$, because it is
hyperbolic and there are no source or sink ($\text{div }X=0$). Also, using that
$\sum \lambda_i =0$ we have that:
\begin{itemize}
\item $\lambda_3(\sigma)<0$ $\Rightarrow$
$-\lambda_3(\sigma)<\lambda_1(\sigma)$,
\item $\lambda_3(\sigma)>0$ $\Rightarrow$
$-\lambda_3(\sigma)>\lambda_2(\sigma)$.
\end{itemize} 
\end{proof}

Now we give a sufficient condition which guarantees that $\Lambda$ is the whole
manifold.

\begin{theorem}
\label{t.nosing}
If $\Lambda$ is a transitive isolated set of a conservative vector field $X$,
such that every $x\in Crit_X(\Lambda)$ is hyperbolic and contains
robustly the unstable manifold of a critical element then $\Lambda=M$.
\end{theorem}

\begin{proof}
We will use the following abstract lemma, which can be found in~\cite{MPP}:

\begin{lemma}[2.7 of~\cite{MPP}]
\label{l.27}
If $\Lambda$ is an isolated set of $X$ with $U$ an isolating block and a
neighborhood $W$ of $\Lambda$ such that $X_t(W)\subset U$ for every $t\geq 0$
then $\Lambda$ is an attracting set of $X$.
\end{lemma}

We will prove that the hypothesis of the previous lemma are satisfied. This will
imply that $\Lambda=M$ since there are no attractors for conservative systems.
If there are no such $W$, then there exists $x_n\to x\in\Lambda$ and $t_n>0$
such that $X_{t_n}\in M-U$ and $X_{t_n}\to q\in \overline{M-U}$. Now, get
$V\subset \ov{O}\subset int(U)\subset U$ a neighborhood of $\Lambda$. We have
$q\notin \ov{V}$.

There exist a neighborhood $\SU_0\subset \SU$ of $X$ such that
$\Lambda_Y(U)\subset V$ for every $Y\in \SU_0$. And by hypothesis we have
$W^u_Y(x_0(Y))\subset V$.

If $x\notin Crit(X)$, then let $z$ such that $\Lambda=\omega_X(z)$, $p\in
W^u_X(x_0)-O_X(x_0)$, so $p\in \Lambda$. Now we have $z_n\in O_X(z)$, $t_n'>0$
such that $z_n\to p$ and $X_{t_n'}(x_n)\to x$. Now, by connecting
lemma~\cite{WX} there exist a $Z\in\SU_0$ such that $q\in W^u_Z(x_0(Z))$ and
$q\notin \ov{V}$ this contradicts the non existence of $W$.

If $x\in Crit(X)$ we can use the Hartman-Grobman theorem to find $x_n'$ in the
positive orbit of $x_n$ ,$r\in (W^u_X(x)-O_X(x))-Crit(X)$ and $t_n'>0$ such that
$x_n'\to r$ and $X_{t_n'}(x_n')\to q$. This is a reduction to the first case
since the following lemma says that $r\in \Lambda$  (the proof of the lemma
below is the same as in conservative case,
since it only use connecting lemma and $\lambda$-lemma):

\begin{lemma}[2.8 of~\cite{MPP}]
\label{l.28}
If $\Lambda$ satisfy the hypothesis of the theorem~\ref{t.nosing} then
$W^u_X(x)\subset \Lambda$ for any $x\in Crit_X(\Lambda)$.
\end{lemma}
\end{proof}

Now we prove that there are no singularities. If $\Lambda=M$ then by
theorem~\ref{t.doering} $X$ will be Anosov, hence whitout singularities. So we
can assume that $\Lambda$ is a proper set with a singularity and find a
contradiction. Now we follow the argument on~\cite{MPP}.

So taking $X$ or $-X$ there exist a singularity $\sigma$ such that
$\dim(W^u_X(\sigma))=1$. Let $U$ an
isolating block such that $\Lambda_Y(U)$ is a connected transitive set for any
$Y\in\SU$. We will prove that $W^u_Y(\sigma(Y))\subset U$ for any $Y\in\SU$, so
by theorem~\ref{t.nosing} we have a contradiction.

Supose that this don't happen for an $Y\in\SU$. By the dimensional hypothesis we
have two branches $w^+$ and $w^-$ of $W^u_X(\sigma)-\{\sigma\}$. Let $q^+\in
w^+$ and $q^-\in w^-$ (also $q^{\pm}(Y)$ the continuation). Because the
negative orbit of $q^{\pm}(Y)$ converge to $\sigma(Y)$ and the unstable
manifold is not containing on $U$ then there exist a $t>0$ such that or
$Y_t(q^+(Y))$ or $Y_t(q^-(Y))$ is not in $U$. We can assume the first case. We
know that there exist a neighborhood $\SU'\subset SU$ of $Y$ such that for any
$Z\in \SU'$ we have $W_t(q^+(Z))\notin U$. Let $z\in \Lambda_Y(U)$ with dense
orbit on $\Lambda_Y(U)$, this implies that $q^-(Y)\in \omega_Y(z)$. So we can
find a sequence $z_n\to q^-(Y)$ in $O_Y(z)$ and $t_n>0$ such that
$Y_{t_n}(z_n)\to q$ for some $q\in W^s_Y(\sigma(Y))-\{\sigma(Y)\}$. We set
$p=q^-(Y)$.

By the connecting lemma, there exist a $Z\in\SU'$ such that $q^-(Z)\in
W^u_Z(\sigma(Z))$, using the same arguments on lemma~\ref{l.28} we can find a
$Z'\in\SU'$ and $t'>0$ such that $Z'_{t'}(q^-(Z'))\notin U$. This shows that
$\sigma(Z')$ is isolated in $\Lambda_{Z'}(U)$ by connecteness we gain that
$\Lambda_{Z'}(U)$ is trivial. A contradiction.

\vspace{.2cm}
Using the same proof of theorem~\ref{t.morales} we can obtain the same result
for conservative vector fields but in dimension greater than 3. So we have that
generically if $X$ has a finite number of homoclinic classes then it is
transitive. Indeed, by Pugh's lemma we obtain $M=Sing(X)\cup
\bigcup\limits_{i=1}^n H(p_i)$ so there are no singularities and only one homoclinic
class, hence transitive. This result cannot be directly strenghted to obtain transitivity of
generic conservative flows, because the connecting lemma of Bonatti,
Crovisier~\cite{BC} is not available in this case. Furthermore, we believe that far away from elliptic points it can be
proved that there is a dominated splitting (the weak Herman's conjecture for
flows). This will need a version for flows of the
Bonatti-Diaz-Pujals~\cite{BDP} theorem.

\section{Robust transitive conservative diffeomorphism}\label{s.dom}

Now we study the existence of a dominated splitting for conservative
diffeomorphisms. This is a topic of many researches and several results on the
dissipative case are well know and we see that transitive systems plays an
important role, for example we have the following result of~\cite{BDP}: \emph{every $C^1$
robustly transitive diffeomorphism has a dominated splitting}. This theorem is
preceding by several results in particular cases : 

\begin{itemize}
\item Ma\~ne~\cite{Ma}: \emph{Every $C^1$-robustly transitive diffeomorphism on
a compact surface is an Anosov system, i.e. it has a hyperbolic (hence dominated) splitting.}
\item Diaz, Pujals, Ures~\cite{DPU}: \emph{There is an open, dense set of 3-dimensional $C^1$-robustly transitive diffeomorphisms admitting dominated splitting}.
\end{itemize}

We use the perturbation lemmas to prove
the conservative version of the theorem by Bonatti, Diaz and Pujals~\cite{BDP} mentioned above:

\begin{maintheorem}\label{BDP conserv.}Let $f$ be a $C^1$-robustly transitive 
conservative diffeomorphism. Then $f$ admits a nontrivial dominated splitting
defined on  the whole $M$.   
\end{maintheorem}

In particular, in dimension 2, we extend Ma\~n\'e's theorem for
conservative robustly transitive systems. 

\begin{theorem}\label{Mane} A $C^1$-robustly transitive conservative system
$f\in\mundo$ is an Anosov diffeomorphism.
\end{theorem}

As remarked in~\cite{BDP}, it is sufficient to show that a conservative robustly
transitive diffeomorphism $f$ can not have periodic points $x\in Per^n(f)$ whose
derivative $D_x f^n$ is the identity. In fact, this is a consequence of :

\begin{proposition}[Corollary 0.4 in~\cite{BDP}]\label{dd de rob transit} Let 
$f\in \mundo$ be a
conservative \emph{transitive} diffeomorphism such that there is a neighborhood 
$\SU\ni f$ such that for any
$g\in\SU$ and every $x\in Per^n(g)$ a periodic point of $g$, the matrix $D_x
g^n$ has at least one eigenvalue of modulus different from one. Then $f$ admits
a dominated splitting. 
\end{proposition}

Also we remark that an immediate corollary of
theorem~\ref{BDP conserv.} is :

\begin{corollary}\label{Ali} Let $f$ be a $C^1$-stably ergodic diffeomorphism in
$\submundo$. Then $f$ admits a dominated splitting.
\end{corollary}

This corollary positively answers a question posed by Tahzibi in~\cite{T}.

\subsection{Proof of Theorem~\ref{BDP conserv.}}

First we give some definitions. We denote by $Per^n(f)$ the periodic points of $f$ with period $n$ and 
$Per^{1\leq n}(f)$ the set of periodic points of $f$ with period at most $n$. 
If $p$ is a hyperbolic saddle of $f$, the homoclinic class of $p$, $H(p,f)$, is the
closure of the transverse intersections of the invariant manifolds of $p$.

Using the $C^2$ pasting lemma~\ref{l.C2pasting} we have immediately the following:

\begin{lemma}
\label{Xia}If $f$ is a $C^2$-conservative diffeomorphism and $x\in
Per^n(f)$ is a periodic point of $f$ such that $D_x f^n$ is the identity matrix
then for any $\ep>0$, there is a $\ep$-$C^{1}$-perturbation $g$ of $f$ 
and some neighborhood $U$ of $x$ satisfying $g^n|_U = id$.      
\end{lemma}

We now prove the theorem~\ref{BDP conserv.}: 

\begin{proof}
If $f$ is robustly transitive then
there is some neighborhood $\SU\ni f$ such that for any $g\in\SU$, $g$ is
transitive. In particular, for any $g\in\SU$, there are not periodic points
$x\in Per^n(g)$ such that $D_x g^n$ is the identity. Indeed, the existence of
such a $g\in\SU$ would imply, by lemma~\ref{Xia}, there is a perturbation
$\widetilde{g}\in\SU$ of $g$ satisfying $(*)\ \widetilde{g}^n|_U=id$, where $U$ is a
small neighborhood of $x$. However, it is easy that $(*)$ is a contradiction
with the transitivity of $\widetilde{g}$. Thus, the proposition~\ref{dd de rob
transit} says $f$ admits a dominated splitting. 
\end{proof}


\section{A conservative version of a result by Carballo-Morales-Pac\'{\i}fico}
\label{cmp}

In this section we deal with a well known open question stated by Michael Herman in the 1998 ICM~\cite{Her}: \emph{A conservative system away from elliptic points has a dominated splitting?}

We deal with a weaker version of this conjecture: \emph{Generically a
conservative system far away from elliptic points has a dominated splitting?}. In dimension 2, this question has an affirmative answer given by Newhouse~\cite{N}: \emph{Generically a $C^1$-conservative diffeomorphism $f$ is Anosov or has a dense set of elliptic points}.

In higher dimensions, we show a conservative version of some recent results
of Carballo, Morales and Pac\'{\i}fico~\cite{CMP} to show this conjecture
with an extra hypothesis: \emph{Generically a $C^2$-conservative diffeomorphism $f$
far away form elliptic points and with a finite number of homoclinic classes has a
dominated splitting}. However, \emph{we want to stress that} Herman's conjecture
was recently proved by
Bonatti-Crovisier, thus this corollary is not a main result of this paper. In fact, as pointed out by ourselves, their connecting lemma
for pseudo-orbits implies the Herman's conjecture in its full version, so our
``corollary'' is nothing but an example of application of the conservative
version of~\cite{CMP}. For more details,
see~\cite{BC}.

\begin{theorem}\label{Herman generico}There is a residual set 
$\SR\subset \mundo$ such
that if $f\in\SR$ is a diffeomorphism with finitely many homoclinic classes and
there is a neighborhood $\SU\ni f$ in $\mundo$ such that for any
$g\in\SU$ and every $x\in Per^n(g)$ a periodic point of $g$, the matrix $D_x
g^n$ has at least one eigenvalue of modulus different from one, then $f$ admits
a dominated splitting and it is a transitive diffeomorphism.
\end{theorem}

\subsection{Proof of theorem~\ref{Herman generico} and a conservative version of
Carballo-Morales-Pacifico}
\label{moraleslemma}

In this subsection, we prove that different homoclinic classes for generic
diffeomorphisms on $\mundo$ are disjoint. We follow the methods used
by~\cite{CMP}, where they prove the same result for generic diffeomorphisms on
$\operatorname{Diff}_{\om}^1(M)$.

The main result of this section is:
\begin{theorem}
\label{t.morales}
For a residual subset $\SF$ of $\mundo$ (conservative or symplectic, 
$\dim(M)\geq 3$), the homoclinic classes of $f\in \SF$ are maximal transitive sets of $f$, therefore different homoclinic classes of $f$ are disjoint.
\end{theorem}

We call a compact subset $\Lp \subset M$ a {\it Lyapunov stable set} for $f$ if every open neighborhood $U\supset \Lp$, there exists an smaller neighborhood $V\supset \Lp$ such that $f^n(V)\subset U$ for all $n\geq 0$. We say that $\Lm$ is a Lyapunov unstable set if it is a Lyapunov stable set for $f^{-1}$. We call $\La$ a {\it neutral set} if $\La=\Lp\cap\Lm$ for some Lyapunov stable set $\Lp$ and Lyapunov unstable set $\Lm$.

We will use the following lemma, which can be found in \cite{CMP} (lemma 2.2)

\begin{lemma}
\label{l.neutral}
If $\Lambda$ is neutral for $f$, then $\Lambda$ is $f$-transitive if and only if $\Lambda$ is a maximal transitive set for $f$. In particular, distinct transitive neutral sets are disjoint.
\end{lemma}


Following~\cite{CMP} we want to prove: 
\begin{theorem}
\label{t.neutral}
There exists a residual set $\SR$ of $\mundo$ where every homoclinic class of any $f\in \SR$ is a neutral set.
\end{theorem}

This together the lemma~\ref{l.neutral} finish the proof of the theorem~\ref{t.morales}.

We follow the steps of \cite{CMP}. For this we need two propositions:

\begin{proposition}
\label{p.critlyap}
There exists a residual set $\SR_1\in \mundo$ where every $f\in \SR_1$ and $p\in Per(f)$, $\overline{W^u(p)}$ (resp. $\overline{W^s(p)}$) is a Lyapunov stable set of $f$ (resp. $f^{-1}$).
\end{proposition}

\begin{proposition}
\label{p.critfecho}
There exists a residual set $\SR_2 \in \mundo$ where every $f\in \SR_2$ satisfies $H_f(p)=\overline{W^u(p)}\cap\overline{W^s(p)}$, for all $p\in Per(f)$.
\end{proposition}

The proof of Theorem~\ref{t.neutral} follows by taken $\SR=\SR_1 \cap \SR_2$ which is also a residual set.

We will need the connecting lemma which can be found in \cite{WX}.

\begin{lemma}
\label{l.xia}
Let $x$ a non periodic point of $f$. For any $C^1_{\omega}$ neighborhood $\SU$ of $f$, there are $\rho>1$, $L\in \N$ and $\eps_0>0$ such that for any $0<\eps <\eps_0$ and any two points $q,t \notin \Delta=\bigcup_{n=1}^{L}f^{-n}(B(x,\delta))$, if exists integers $m_1>0$ and $m_2< 0$ such that 
$f^{m_1}(q)$ and $f^{m_2}(t)$ are inside of $B(x,\delta/\rho)$ then there exists $g\in \SU$ such that $g=f$ on $\Delta^c$ and $t$ is on the positive $g$-orbit of $q$.
\end{lemma}

We will fix some notations. We denote $\com$ the space of all compact subsets of $M$ with the Hausdorff Topology.
Also we denote by $\SK$ the set of all $f \in \mundo$ which are Kupka-Smale. This set is a generic set by a theorem of Robinson~\cite{Rob}.

\smallskip

{\bf 1. Proof of Lemma~\ref{p.critlyap}}

First we prove the local version of the lemma~\ref{p.critlyap}. Then, it follows from the same methods used on the proof of lemma 3.5 of~\cite{CMP}.

\begin{lemma}
\label{l.loccritlyap}
Let $f\in \SK$ and $n\leq 0$, then there exists a neighborhood $\SU(f,n)\in\mundo$ of $f$ and a residual set $\SR(f,n)\subset \SU(f,n)$ such that for every $g\in \SR(f,n)$, any $g$-periodic point $p$ of period $m\leq n$ the set $\overline{W^u(p)}$ (resp. $\overline{W^s(p)}$)  is Lyapunov stable for $g$ (resp. $g^{-1}$).
\end{lemma}

\begin{proof}
The proof has the same lines of Lemma 5.1 from \cite{CMP}, so we sketched the proof. If $\SU(f,n)$ is small, by hyperbolicity we have that $Per^{1\leq n}(g)=\{p_1(g),...,p_k(g)\}$ for any $g\in \SU(f,n)$. For $i=1,..k$ we set the function $\Psi_i:\SU(f,n) \rightarrow \com$ by: 
$$\Psi_i(g)=\overline{W^u_g(p_i(g))}.$$
This map is lower semi-continuous on the Hausdorff topology. Then is continuous on a residual set $\SR_i\in\SU(f,n)$, define $\SR=\SK\cap(\cap_i\SR_i)$.

Now, get $g\in \SR$ and $p=p_i(g)\in Per^m(g)$, then we claim that$\Psi_i(g)$ is Lyapunov stable for $g$. 

If not, by definition there exists an open set $U\supset \Psi_i(g)$, $x\in \Psi_i(g)$ and sequences $x_i\to x$, $n_i\geq0$ such that $f^{n_i}(x_i)\notin U$. Now we have two cases, if $x$ is a periodic point or not.

By an argument of \cite[p.~14]{CMP} is enough show a contradiction in the case of $x$ is not a periodic point, in fact if $x$ is a periodic point, using the fact of $g$ is Kupka-Smale and Hartman-Grobman theorem, we can found another sequences $x_i'\to x'$ and $n_i'\geq 0$ like above, with $x'\in W^u_g(x)$, $x'$ is not on the orbit of $x$ and $g^{n_i'}(x_i')\notin U$.

Now, let's focus on the case of $x\notin Per(g)$. Because $\Psi_i$ is continuous on g, and $\overline{\Sun}\subset U$, there is a $C^1_{\omega}$-neighborhood $V$ of $g$ inside of $\SU(f,n)$, such that $\overline{W^u_h(p_h)}\subset U$, for any $h$ in $V$, where $p_h$ is the continuation of $p$. Now we use the lemma~\ref{l.xia} for $x$, $g$ and $\SU$ giving us $\rho$, $L$ and $\eps_0$.

By the non-periodicity of $x$, $f^{-i}(x)$ do not intersect the orbit of $p$ for $i=0,..,L$, furthermore, there exists an $\eps \leq \eps_0$ such that $\Delta=\bigcup_{n=1}^{L}f^{-n}(B(x,\eps))$ do not intersect the orbit of $p$, and $\Delta \subset U$. Now choose an open set $S$ containing the orbit of $p$ such that: $\overline{S}\subset U$ and $V\cap \Delta = \emptyset$ (*).

Let $i$ (large enough) such that $x_i\in B(x,\eps \backslash \rho)$ and set  $t=g^{n_i}(x_i) \notin U$. Then exists an $m_2< 0$ such that $f^{m_2}(t)\in B(x,\eps \backslash \rho)$. Now get $q\in (W^u(p)\backslash\{p\})\cap S$ such that there exists a $m_1>0$ which $f^{m_1}(q)\in B(x,\eps\backslash\rho)$ and $f^{-i}(q) \in S$ for any $i>0$ (therefore by (*) we have $f^{-i}(q) \notin \Delta$ (**)). 

Now we apply the lemma~\ref{l.xia} and find an $h\in V$, $h$ equal to $g$ in
$\Delta^c$. Then by (*), (**), we have $p_h=p$ and $q\in W^u_h(p)$. But the lemma also says that $t$ is in the $h$-positive orbit of $q$, then $t\in W^u_h(p)$ and $t\notin U$, and this is a contradiction with the fact that for a $h\in V$, we have $\overline {W^u_h(p_h)}\subset U$. 
\end{proof}

The proof of proposition~\ref{p.critlyap} follows from the lemma above by standard methods (see the proof of lemma
3.5, p. 12 on \cite{CMP}).

\smallskip

{\bf 2. Proof of Lemma~\ref{p.critfecho}}

Again the proof follows from a local form of the proposition and the methods of the proof of lemma 3.5 in~\cite{CMP}.

\begin{lemma}
\label{l. localfecho}
Let $f\in \SK$ and $n>0$, there exists a neighborhood $\SV(f,n)$ of $f$ and a residual ser $\SP(f,n)\subset \SV(f,n)$ such
that if $g\in \SP(f,n)$ and $p\in Per^m(g)$ with $m\leq n$ then $H_g(p)=\overline{W^u_g(p)}\cap\overline{W^s_g(p)}$.
\end{lemma}

\begin{proof}
Let $\SV(f,n)$ neighborhood of $f$ such that $Per^{1\leq n}(g)=\{p_1(g),...,p_k(g)\}$ on $\SV(f,n)$, and again let
$\Psi_i(g)=H_g(p_i(g))$ for $i=1,..k$, a lower semicontinuous function, and $\SP^i(f,n)$ the set of points of continuity of
$\Psi_i$ on $\SV(f,n)$ (a residual set). Set $\SP=\SK\cap(\cap_i\SP^i(f,n))\cap\SR$ (residual) where $\SR_1$ is the residual set
given by proposition~\ref{p.critlyap}.

Let $p=p_i(g)\in Per^m(g)$ for some $i$ with $m\leq n$ and $g\in \SP(f,n)$. Suppose that the lemma is false for $\Psi_i(g)$, then there is $x\in \overline{W^u_g(p)}\cap\overline{W^s_g(p)}\backslash H_g(p)$.

Again we have two cases, $x$ be a periodic point or not, like in lemma~\ref{l.loccritlyap} we can suppose $x$ a non periodic point. We fix then $K$ a compact neighborhood of $x$ such that $K\cap H_g(p)=\emptyset$, then by continuity there is $\SU$ neighborhood of $g$ where $K\cap  H_h(p_h)=\emptyset$ for all $h\in \SU$.

Let $\rho$, $L$ and $\eps_0$ in the lemma~\ref{l.xia}. We found an $\eps<\eps_0$ such that $g^i(p)\notin \Delta=\bigcup_{i=0}^{L}f^{-i}B(x,\eps)$ and $B(x,\eps)\subset K$. Let $V$ containing the orbit of $p$ and $V\cap \Delta=\emptyset$. Again we can find $q\in W^u_g(p)\backslash \{p\}\cap V$, $t\in W^s_g(p)\backslash \{p\}\cap V$, $m_1>0$ and $m_2<0$ like in  lemma~\ref{l.xia} and also, with $g^{-i}(q)\in V$ and $g^i(t)\in V$ for all $i>0$. 

If there exists $i>0$ with $g^i(q)=t$ then q be an homoclinic point of $g$ passing on $K$ absurd, so there is not such $i$.

So by the lemma~\ref{l.xia} there exists $h\in \SU$ equal to $g$ off $\Delta$ an
$t$ is in the future orbit of $q$. Like before, we have $p_h=p$, $q\in W^u_h(p)$
and $t\in W^s_h(p)$, by a small perturbation we can assume that the $h$-orbit of
$q$ is in $H_h(p_h)$. This implies $K\cap H_h(p_h)\neq \emptyset$, a contradiction. 
\end{proof}

Now we prove the theorem~\ref{Herman generico} :

\begin{proof}
Using lemma 7.8 of~\cite{BDP}, we
know that non-trivial homoclinic classes of periodic hyperbolic points are dense
in $M$, for generic diffeomorphisms $f\in \mundo$. Moreover, the
theorem 6 in~\cite{BDP} and the assumption that there is not periodic point
with all its eigenvalues of the derivative equal to one, for any 
diffeomorphism $g$ $C^1$-close to $f$, imply that $M$ is the
union of the invariant compact sets $\La_i(g)$ defined by the closure of the
union of the homoclinic classes with a dominated decomposition of index $i$.
But, the theorem~\ref{t.morales} implies that for generic
diffeomorphisms, different homoclinic classes are disjoint. In particular, by
the finiteness of homoclinic classes, $M$ is the \emph{disjoint} union of a 
finite number of homoclinic classes (which are invariant compact sets with a
dominated splitting). Because $M$ is connected, $M$ is exactly one homoclinic
class. Since, in our case, homoclinic classes are transitive sets and 
admit a dominated splitting, $f$ admits a dominated splitting and $f$ is
transitive. 
\end{proof}


\bibliographystyle{alpha}
\bibliography{bib}

\begin{thebibliography}{35}

\bibitem[Au]{Au}
T.~Aubin,
\newblock Nonlinear Analysis on Manifolds. Monge-Amp\`ere Equations,
\newblock Springer-Verlag, 1982. 

\bibitem[BFP]{BFP}
J.~Bochi, B.~Fayad and E.~Pujals.
\newblock \ \ In preparation.

\bibitem[BC]{BC}
C. Bonatti and S. Crovisier,
\newblock R\'ecurrence et G\'en\'ericit\'e,
\newblock Preprint 2003.

\bibitem[BDP]{BDP}
C.~Bonatti, L.~Diaz, E.~Pujals,
\newblock A $C^1$-generic dichotomy for diffeomorphisms: Weak forms of
hyperbolicity or infinitely many sinks or sources,
\newblock Preprint 2002, to appear in \emph{Annals of Math}.

\bibitem[BS]{BS}
V. Biragov and L. Shilnikov,
\newblock On the bifurcation of a saddle-focus separaratrix loop in a
three-Dimensional conservative dynamical system,
\newblock \emph{Selecta Mathematica Sovietica}, v.11, no.4, 333--340, 1992.

\bibitem[BV]{Bochi-Viana}
J.~Bochi and M.~Viana,
\newblock The Lyapounov exponents of generic volume preserving and symplectic
systems, 
\newblock Preprint 2002, to appear in \emph{Annals of Math.}

\bibitem[CMP]{CMP}
C.~Carballo, C.~Morales and M.~Pac\'{i}fico,
\newblock Homoclinic classes for generic $C^1$ vector fields, 
\newblock Preprint 2001, to appear in \emph{Erg. Theory and Dyn. Systems}.

\bibitem[DM]{DM}
B.~Dacorogna, J.~Moser,
\newblock On a partial differential equation involving the jacobian determinant,
\newblock \emph{Ann. Inst. Poincar\'e}, V.7. 1--26. 1990. 

\bibitem[DPU]{DPU}
L.~Diaz, E.~Pujals, R~Ures,
\newblock Partial hyperbolic and robust transitivity.
\newblock \emph{Acta. Math.}, v. 183, 1--42, 1999.

\bibitem[D]{D}
C.~Doering.
\newblock Persistently transitive vector fields on three-dimensional manifolds,
\newblock \emph{Proc. on Dynamical Systems and Bifurcation Theory}, Pitman Res.
Notes Math. Ser., 160, 59--89, 1987.

\bibitem[GT]{GT}
D.~Gilbarg and N.~Trudinger
\newblock Elliptic Partial Differential Equations of Second Order,
\newblock Springer-Verlag, 1998.

\bibitem[Her]{Her}
M.~Herman, 
\newblock Some open problems in dynamical systems.
\newblock \emph{Doc. Mathematica. Extra Volume ICM}, 1998-II, 797--808.

\bibitem[H]{H}
L.~H\"ormander,
\newblock Linear Partial Differential Operators,
\newblock Springer-Verlag, 1976. 

\bibitem[LU]{LU}
O.~Ladyzenskaya and N.~Uralsteva,
\newblock Linear and Quasilinear Elliptic Partial Differential Equations,
\newblock Academic Press, 1968. 

\bibitem[Ma]{Ma}
R.~Ma\~n\'e,
\newblock An ergodic closing lemma.
\newblock \emph{Ann. Math.}, v. 116, 503--540. 1982.

\bibitem[MPP]{MPP}
C.A. Morales, M.J. Pac\'{\i}fico, E. Pujals.
\newblock Robust transitive singular sets for 3-flows are partially hyperbolic
attractors or repellers.
\newblock preprint. 2002.

\bibitem[Mo]{M}
J.~Moser,
\newblock On the volume element of a manifold.
\newblock \emph{Trans. AMS}, v.120, 286--294, 1965.

\bibitem[N]{N}
S.~Newhouse,
\newblock Quasi-elliptic periodic points in conservative dynamical systems.
\newblock \emph{Amer. J. Math.}, V. 99, no. 5, 1061--1087, 1977.

\bibitem[R]{Rob}
R.C.~Robinson,
\newblock Generic Properties of conservative systems, I and II,
\newblock \emph{Amer. J. Math}, 92, 562--603 and 897--906, 1970.

\bibitem[RY]{RY}
T.~Riviere and D.~Ye, 
\newblock Resolutions of the prescribed volume form equation,
\newblock \emph{Nonlinear Diff. Eq. Appl.}, 3, 323--369, 1996.

\bibitem[T]{T}
A.~Tahzibi,
\newblock Stably ergodic systems which are not partially hyperbolic,
\newblock Preprint 2002.

\bibitem[V]{Max}
M.~Viana,
\newblock What's new on Lorenz Strange Attractors
\newblock \emph{Math. Inteligencer}, vol. 22, no. 3, 6--19, 2000. 

\bibitem[WX]{WX}
L.~Wen and Z.~Xia,
\newblock $C^1$ connecting lemmas,
\newblock \emph{Trans. Amer. Math. Soc.}, 352, no. 11, 5213--5230, 2000.

\bibitem[X]{X}
Z.~Xia,
\newblock Homoclinic points in symplectic and volume preserving diffeomorphisms,
\newblock \emph{Comm. in Math. Phys.}, Vol.177 , 435--449, 1996.

\bibitem[Z]{Z}
E.~Zehnder,
\newblock Note on smoothing symplectic and volume preserving diffeomorphisms,
\newblock \emph{Lect. Notes in Math.}, 597, 828--854, 1977.

\bibitem[Zu]{Zu}
C.~Zuppa,
\newblock Regularisation $C^{\infty}$ des champs vectoriels qui pr\'eservent
l'el\'ement de volume.
\newblock \emph{Bol. Soc. Brasileira Matem.}, 10 (2), 51--56, 1979.

\end{thebibliography}

\vspace{1cm}

\noindent 
		\textbf{Alexander Arbieto} ( alexande{\@@}impa.br )\\
		\textbf{Carlos Matheus} ( matheus{\@@}impa.br )\\
		IMPA, Est. D. Castorina, 110, Jardim Bot\^anico, 22460-320\\
		Rio de Janeiro, RJ, Brazil

\end{document}